\documentclass[12pt]{amsart}
\usepackage{amscd, amssymb}
 
\newtheorem{pr}{Proposition}
 
\newtheorem{lm}{Lemma}
\newtheorem{tm}{Theorem}
\newtheorem{cor}{Corollary}

\newcommand{\proj}{\mathbb P}

\newcommand{\barr}{\overline}
\newcommand{\rarr}{\rightarrow}

\newcommand{\com}{\mathbb{C}}

\newcommand{\Z}{\mathbb{Z}}

\newcommand{\Lie }{{\mathrm{Lie}}}
\newcommand{\Tan}{{\mathrm{Tan}}}

\newcommand{\GG}{\mathbf{G}}
\newcommand{\PP}{\mathbf{P}}
\newcommand{\TT}{\mathbf{T}}
\newcommand{\BB}{\mathbf{B}}
\newcommand{\WW}{\mathbf{W}}

\newcommand{\eqq}{\stackrel{\sim}{=}}

\newcommand{\deli}{\bigtriangleup}

\newcommand{\bpf}{\noindent {\em Proof.} }
\newcommand{\epf}{\qed \vspace{+10pt}}

\begin{document}

\title
{The connectedness of the moduli space of maps to homogeneous spaces}

\author{B. Kim}
\address{Pohang University of Science and Technology}
\email{bumsig@postech.ac.kr}
\author{R. Pandharipande}
\address{California Institute of Technology}
\email{rahulp@cco.caltech.edu}
\date{25 March 2000}
\maketitle

\pagestyle{plain}
\setcounter{section}{-1}
\section{Introduction}
\label{yintro}
Let $X$ be a compact algebraic homogeneous space: $X=\GG/\PP$
where 
$\GG$ is a connected complex semisimple algebraic group
and $\PP$ is a parabolic subgroup. Let $\beta\in H_2(X, \mathbb{Z})$.
The (coarse) moduli space $\overline{M}_{g,n}(X,\beta)$
of $n$-pointed genus $g$ stable maps parameterizes 
the data
$$[\mu: C\rarr X, p_1, \ldots, p_n ]$$  satisfying:
\begin{enumerate}
\item[(i)]
$C$ is a complex, projective, connected, reduced, (at worst) nodal curve
of arithmetic genus $g$.
\item[(ii)] The points $p_i \in C$ are distinct and lie in 
           the nonsingular locus.
\item[(iii)] $\mu_*[C]=\beta$.
\item[(iv)]  The pointed map $\mu$ has no infinitesimal automorphisms.
\end{enumerate}
Since $X$ is convex, the genus 0  moduli
space $\overline{M}_{0,n}
(X, \beta)$ 
is of pure dimension $$\text{dim}(X)+ \int_\beta c_1(T_X)+ n-3.$$
Moreover, $\overline{M}_{0,n}
(X, \beta)$
is locally
the quotient of a nonsingular variety by a finite
group. 
For general $g$, the space $\overline{M}_{g,n}(X, \beta)$
may have singular components of different dimensions.
Stable maps in algebraic geometry were first defined in [Ko].
Basic properties of the moduli space $\overline{M}_{g,n}(X,\beta)$
can be found in [BM], [FP], and  [KoM].
The
following connectedness result is proven here.
\begin{tm} 
$\overline{M}_{g,n}(\GG/\PP, \beta)$ is a connected variety.
\end{tm}

\noindent
This result may be viewed as analogous to the connectedness
of the Hilbert scheme of projective space proven by Hartshorne.
As in [Har], connectedness is obtained via maximal degenerations.

Since
$\overline{M}_{0,n}(X, \beta)$ 
has quotient singularities, connectedness is equivalent
to irreducibility.
\begin{cor}
\label{fred}
$\overline{M}_{0,n}(\GG/\PP, \beta)$ is an irreducible variety.
\end{cor}

Corollary \ref{fred} is easy to verify
in case $X$ is a projective space. When $X$ is a Grassmannian,
the irreducibility follows from Str\o mme's Quot scheme
analysis [S].
A proof of Corollary \ref{fred} can be
found in case $\GG=\mathbf{SL}$ in [MM].
For the variety of partial 
flags in $\com^n$, a proof of irreducibility 
using flag-Quot schemes is established in [Ki].
Results of Harder closely related to Corollary \ref{fred}
appear in [Ha].  
There is an independent proof by J. Thomsen for 
the irreducibility of $\overline{M}_{0,n}(\GG/\PP ,\beta )$
in [T].

\bigskip

The moduli space $\overline{M}_{g,n}(X,\beta)$ has a natural locally 
closed decomposition
indexed by 
stable, pointed, modular graphs $\tau$ (see [BM]).
The strata correspond to maps with domain curves of a fixed topological
type and a fixed 
distribution $\beta_\tau$ of $\beta$. The graph $\tau$ determines
a complete 
moduli space of stable maps $$\overline{M}_{\tau,n}(X,\beta_\tau)$$
together with a canonical morphism:
\begin{equation}
\label{wwww}
\pi_\tau: 
\overline{M}_{\tau,n}(X,\beta_\tau) \rarr \overline{M}_{g,n}(X,\beta).
\end{equation}
A closed decomposition is determined by the images of
these morphisms (\ref{wwww}).
Theorem 1 is a special case of the
following result.

\begin{tm}
$\barr{M}_{\tau,n}(\GG/\PP, \beta_{\tau})$ is a connected variety.
\end{tm}
\noindent
Since $\barr{M}_{\tau,n}(X, \beta_{\tau})$ is
normal in the genus 0 case, we obtain the corresponding corollary.

\begin{cor}
Let $g=0$.
$\barr{M}_{{\tau}, n}(\GG/\PP, \beta_{\tau})$ is an 
irreducible  variety.
\end{cor}
\noindent In particular, all the boundary divisors
of $\overline{M}_{0,n}(X, \beta)$ are irreducible.

Theorem 2
is proven  by studying the maximal torus action on $X$.
The method is to degenerate a general $\GG$-translate
of a map $\mu: C \rarr X$ onto a canonical
$1$-dimensional configuration of $\proj^1$'s in $X$
determined by the 
maximal torus and the Bialynicki-Birula stratification of $X$.

\bigskip

In the genus 0 case, we study 
the Bialynicki-Birula stratification of $\overline{M}_{0,n}(X,\beta)$.
The following result is then deduced from the
rationality of torus fixed components.

\begin{tm}
$\overline{M}_{0,n}(\GG/\PP , \beta )$ is 
rational.
\end{tm}

\noindent
The fixed component rationality is equivalent to a rationality result for
certain quotients of  $\mathbf{SL}_2$-representations proven by Katsylo
and Bogomolov [Ka], [Bog]. It should be noted that the
fixed components will in general be contained in the boundary of the
moduli space of maps -- the compactifaction by stable maps
therefore plays an important role in the proof.

The rationality of the Hilbert schemes of rational curves in
projective space (birational to $\overline{M}_{0,0}(\proj^r,d)$)
is a consequence of Katsylo's results [Ka] and was also studied
by Hirschowitz in [Hi]. 
\bigskip

The main part of this paper was written in 1996 at 
the Mittag-Leffler Institute where the
authors benefitted from discussions with many members. 
Thanks are especially due to I. Ciocan-Fontanine, B. Fantechi,
W. Fulton, T. Graber, and
B. Totaro. Conversations with F. Bogomolov were also helpful.
B. ~K. ~was partially supported by
KOSEF grant 1999-2-102-003-5 and POSTECH grant 1999.
R.~P.~was partially supported by NSF grant 
DMS-9801574 and
an A.~P.~Sloan foundation
fellowship.

\section{The torus action on $\GG/\PP$ }  
\label{torr}
Let $\GG$ be a connected complex semisimple algebraic group.
Let $\PP$ be a parabolic subgroup. 
Select a maximal algebraic torus $\TT$ and
Borel subgroup $\BB$ of $\GG$ satisfying:
  $$\TT \subset \BB \subset \PP \subset \GG.$$
Let $(\GG/\PP)^\TT$ denote the fixed point set of 
the left $\TT$-action on $\GG/\PP$. Three special properties
of this $\TT$-action
will be needed:
\begin{enumerate}
\item[(i)] The $\TT$-action has isolated fixed points.
\item[(ii)] For every point $p\in (\GG/\PP)^\TT$, there
   exits a $\TT$-invariant open set $U_p$ containing
   $p$ which is $\TT$-equivalent to a vector space
   representation of $\TT$.
\item[(iii)] Let $\com^*\subset \TT$ correspond
to an interior point of a Weyl chamber. Then, $(\GG/\PP)^{\com^*}
=(\GG/\PP)^\TT$,
and  the Bialynicki-Birula decomposition obtained from
 the $\com^*$-action is an affine
 {\em stratification} of $\GG/\PP$.

\end{enumerate}
A stratification is a decomposition such that
the closures of the strata are unions of strata.
In general, the Bialynicki-Birula decomposition 
obtained from a $\com^*$-action on a nonsingular variety
need not be
a stratification.

The claims (i)-(iii) are well known. Only a brief summary
of the arguments will be presented here.
Let $\WW$ be the Weyl group of $\GG$ relative to $\TT$.

\begin{lm}  $|(\GG/\BB) ^\TT| = |\WW|$, and $\WW$ 
acts transitively on $(\GG/\BB)^\TT$.
\end{lm}

\bpf   See, for example, [Bor].
\epf

\noindent In particular, $(\GG/\BB)^\TT$ is a finite set.

\begin{lm} The natural map $(\GG/\BB) ^\TT \rarr (\GG/\PP)^\TT$
           is surjective.
\end{lm}

\bpf    
            Let $p \in (\GG/\PP) ^\TT$. 
           The invariant fiber (isomorphic to $\PP/\BB$) over
          the fixed point $p$ is a nonsingular projective variety, and
            hence contains a $\TT$-fixed point by the Borel
            fixed point theorem (or, alternatively,
           this is a Hamiltonian
            action on a compact manifold). 
\epf

\noindent Therefore, $\WW$ acts transitively on the finite
         set $(\GG/\PP)^\TT$.

 A  representation  $\psi:\TT\rarr \mathbf{GL}(V)$ is
   {\em fully definite} if there
   exists a $\com ^*$-basis of $\TT$ for which all
   the weights of the representation are positive integers.
   Equivalently, a fully definite representation can be
   written 
$$\psi(t_1, \ldots, t_r) v_j= \prod_{i=1}^{r} t_i^{\lambda_{ij}} \cdot v_j$$
 where $\lambda_{ij}> 0$ 
  for some choice of $\com^*$-basis of $\TT$ and $\com$-basis $\{v_j\}$
  of $V$.

The point $1\in \GG/\BB$ corresponding to the identity element of $\GG$
is a $\TT$-fixed point.
   The $\TT$-action induces a representation $$\phi:\TT \rarr
   \mathbf{GL}( \text{Tan}_1 \GG/\BB).$$

\begin{lm} The representation
$\phi$ is fully definite.
\label{ddef}
\end{lm}

\bpf      The natural quotient map $q:\GG \rarr \GG/\BB$
          is $\TT$-equivariant for the conjugation action on $\GG$
          and the left action on $\GG/\BB$.
          The differential of $q$ yields an isomorphism from 
          the Adjoint representation of $\TT$ on $\Lie (\GG)/\Lie (\BB)$
          to $\phi$. 
            $\Lie (\GG)/ \Lie (\BB)$ is the space of positive
            roots. This $\TT$-representation space has
            $n$ simple roots (where $n$ is the
            rank of $\GG$). All the 1-dimensional representations
            in $\Lie(\GG)/\Lie(\BB)$ are non-negative tensor products
            of these simple  roots. Moreover, the
            $n$ weight vectors of these simple roots
            are independent in the lattice of
            1-dimensional representations of the torus $\TT$. 
            Lemma \ref{ddef} now follows from Lemma \ref{bbas} below.
\epf

\begin{lm}      Let $\psi:\TT \rightarrow \mathbf{GL}(\com^n)$ 
               be an $n$ dimensional
                   representation of a rank $n$ torus $\TT$.
                   If the $n \times n$ matrix of weights is
                   nonsingular, then the representation
                   is fully definite.
\label{bbas}
\end{lm}
\bpf See [Bi].
\epf

\begin{lm}   The $\TT$-representation  $\text{\em{Tan}} _1 \GG/\PP$
            is fully definite.
\label{dxx}
\end{lm}

\bpf There is a surjection of $\TT$-modules given
            by the differential $\Tan _1 \GG/\BB \rarr \Tan _1 \GG/\PP$.

\epf

\begin{pr}  For every $p \in (\GG/\PP) ^\TT$, there exists
            a $\TT$-invariant Zariski open set $U_p\subset \GG/\PP$ of $p$
            which is $\TT$-equivalent to a vector space representation
            of $\TT$.
\label{rrr}
\end{pr}
\bpf    By a theorem of Bialynicki-Birula [Bi], it suffices to
            show the tangent representation of $\TT$ is fully definite
            at $p$. This is a consequence of Lemma \ref{dxx}
            and the transitivity of the $\WW$-action on $(\GG/\PP)^\TT$.
            (In fact, only definiteness of the tangent representation
             is needed in [Bi].)
\epf

Let $\com^*\subset \TT$ correspond
to an interior point of a Weyl chamber. 
By the analysis of the tangent representation $\phi$, every point of   
$(\GG/\BB)^\TT$ is an isolated fixed point of $\com^*$. The equality 
 $(\GG/\BB)^{\com^*}=(\GG/\BB)^\TT$ follows. Since the map
$(\GG/\BB)^{\com^*} \rarr (\GG/\PP)^{\com^*}$ is surjective, 
$(\GG/\PP)^{\com^*}=(\GG/\PP)^{\TT}$.

For each $p\in (\GG/\PP)^\TT$, let $A_p$ be the set of points
$x\in \GG/\PP$ such that $$ \lim_{ t\rarr 0} tx =p.$$
By Proposition \ref{rrr}, $A_p$ is isomorphic to
the affine space $\com^{r_p}$ where $r_p$ is
the number of positive weights in the
$\com^*$-representation $\text{Tan}_p \GG/\PP$.
The set $\{A_p\}$ is the Bialynicki-Birula
affine decomposition of $\GG/\PP$.
In fact, $\{A_p\}$ coincides (up to the Weyl group action)
with the (open) Schubert cell stratification of
$\GG/\PP$. This is essentially proven in [Bor] 
for the case $\GG/\BB$. The general case $\GG/\PP$ is proven in
[A]. Therefore, $\{A_p\}$ is a stratification.

\section{The $\com^*$-flow}
\label{cflow}
   Let $\com^* \subset \TT$ correspond
to an interior point of a Weyl chamber. 
   Let $s,x_1, \ldots, x_l\in (\GG/\PP)^\TT$ be the fixed points
   corresponding to the unique maximal dimensional stratum
   $A_s$ and the complete set of
   codimension 1 strata, $A_1, \ldots, A_l$, respectively. 
   The points of $A_s$ flow $(t\rarr 0)$ to $s$, and
   the points of  $A_i$ flow $(t\rarr 0)$ to $x_i$.
   Let $U= A_s\cup A_1 \cup \ldots \cup A_l$. Since the
   Bialynicki-Birula decomposition  $\{A_p\}$ is a stratification, 
   $U$ is a Zariski open set with complement
   of codimension at least 2.

   The inverse action of $\com^*$ on $\GG/\PP$ 
   is also a torus action on $\GG/\PP$
   with the same fixed point set.
   Let $A'_s, A'_1, \ldots, A'_l$ be the affine strata for the
   inverse action
   corresponding to the fixed points $s,x_1, \ldots, x_l$.
Let $\text{dim}(\GG/\PP)=m$.
   Since,
$$\text{dim}(A_p)+ \text{dim}(A'_p)= m,$$
   $A'_1, \cdots, A'_l$ are the complete
   set of $1$-dimensional strata for the inverse action.
   Moreover, the closure $P_i=\barr{A'}_i$ can contain only
   the unique $0$-dimensional stratum $A'_s= s$.
   We have shown the closures $P_i$ are contained in $U$. 
   Each $P_i$ is 
   isomorphic to $\proj^1$ (Chevelley [C] proves
   the closed Schubert cells have singularities in codimension
   at least 2).
   The intersection pairing
$$ P_i \cap \overline{A}_j = \delta(i-j)$$
follows from the above analysis.  
   Since the closed strata of the inverse action
   freely generate the integral homology, the classes
   $$[P_1], \ldots, [P_l] \in H_2(\GG/\PP, \mathbb{Z})$$
   span an integral basis of $H_2(\GG/\PP,\mathbb{Z})$.

Let $f:C\rarr \GG/\PP$ be a 
non-constant stable map satisfying the following properties:
\begin{enumerate} 
\item[(i)] The image $f(C)$ lies in $U$.
\item[(ii)] $C$ intersects (via $f$) the divisors $A_i$
transversely at nonsingular points of $C$. 
\item[(iii)] All the markings of $C$ have image in
$A_s$.
\end{enumerate}

If $[f]$ represents the class
$$\beta = \sum_{i=1}^l a_i [P_i] \in H_2(\GG /\PP,\Z),$$ then let
$C$ meet $A_i$ at the $a_i$ distinct points 
$$\{ x_{i,1}, \ldots, x_{i,a_i}\}.$$

We will study the induced
$\com^*$-action on $\overline{M}_{g,n}(\GG/\PP,\beta)$
by translation of maps.
Let $F: C_0 \rarr \GG/\PP$ be the limit in the space of stable maps,
$$F=\lim _{t\rarr 0} \ tf$$  
where $t \in \com^*$.

Define a map $\tilde{F}: \tilde{C} \rarr \GG/\PP$ as follows.
Let the domain $\tilde{C}$ be:
$$\tilde{C} =
C\ \cup \ \bigcup _{i=1}^{l}\ (\cup _{j=1}^{a_i}\PP _{i,j} ^1)$$ where
$\PP _{i,j}^1$ is 
a projective line
attached to $C$ at the point $x_{i,j}$.
Let the markings of $\tilde{C}$ coincide with the
markings of $C$ (note the markings of $C$
 are disjoint from the set $\{x_{i,j}\}$ by condition (ii)).
Define $\tilde{F}$ by
$\tilde{F}(C\subset \tilde{C})=s$ and 
$$ \tilde{F}|_{\PP ^1_{i,j}}:  \PP^1_{i,j} \eqq P_i$$
for each $i$ and $j$.

\begin{pr}
\label{pr2}
If $f$ satisfies conditions (i-iii), then
the $t\rarr 0$ limit $F$ equals 
the stabilization of $\tilde{F}$.
\end{pr}

\bpf 
Let $\deli^\circ\subset \deli$ be the punctured holomorphic
disk at the origin.
Let $$h: C \times \deli^\circ \rarr \GG/\PP$$ be the map defined
by $h(c,t)=tf(c)$.
The $\com^*$-action on $A_s$ extends to a map
$$\com \times A_s \rarr A_s$$
since the $\com^*$-action on $A_s$ is a vector space representation
with positive weights.
The map $h$ thus extends to a 
map 
$$h:C\times\deli \setminus  \{ x_{i,j}\times 0\} \rarr \GG/\PP$$
since
the $f$-image
of $C\setminus \{ x_{i,j}\}$ lies in $A_s$.
Note,
\begin{equation}
\label{wgw}
h(C \setminus \{ x_{i,j}\}, 0  )=s.
\end{equation}
After a suitable blow-up $$\gamma:S \rarr C\times \deli$$
 supported along the isolated nonsingular 
points $\{x_{i,j} \times 0\}$ 
of $C\times \deli$,
there is a morphism $h': S\rarr \GG/\PP$. 

The limit as $t\rarr 0$ of $tf(x_{i,j})$ equals $x_i$.
Hence,
the exceptional divisor $C_{i,j}$ of $\gamma$ over  $x_{i,j}$ 
connects the points $x_i$ to $s$ under the map $h'$.
The image $h'(C_{i,j})$ 
thus represents an effective curve class containing the class
$[P_i]$.
By degree considerations over all the exceptional
divisors $C_{i,j}$,
we conclude $h'(C_{i,j})$ is of curve class exactly $[P_i]$.
As $P_i$ is the unique $\com^*$-fixed curve of class $[P_i]$ connecting
the points $x_i$ and $s$, 
$$h'(C_{i,j})=P_i.$$ 
We may assume $S$ to be nonsingular (away from
the original nodes of $C$) and each  $C_{i,j}$ to be a 
normal crossings divisor -- possibly
after further blow-ups and base changes altering
only the special fiber over $0\in \deli$.
We then conclude  each $C_{i,j}$ has
a single component which is mapped to $P_i$ isomorphically (and
the other components of $C_{i,j}$ are contracted). 

After blowing-down
the $h'$-contracted components of each $C_{i,j}$,
we obtain a map $h'':S'' \rarr \GG/\PP$
which is a family of nodal maps over $\deli$.
The fiber of $S''$ over $t=0$ is isomorphic to $\tilde{C}$.
Moreover, the condition
$\tilde{F}(C\subset \tilde{C})=s$ follows directly from (\ref{wgw}).

The limit {stable} map $F$ is then simply obtained by
stabilizing the map $\tilde{F}$. We have carried out the
stable reduction of the family of maps $tf$ (see [FP]).
\epf

\section{Connectedness}
\label{conn}
Let $[\mu]$ denote the point
$[\mu:C\rarr X, p_1, \ldots, p_n]\in \overline{M}_{g,n}(X,\beta).$
The stable, pointed, modular graph $\tau$ 
with $H_2(X,\mathbb{Z})$-structure canonically associated to 
$[\mu]$ consists of the following data:
\begin{enumerate}
\item[(i)]
The pointed dual
graph of $C$:
\begin{enumerate}
\item[(a)] The vertices $V_\tau$ 
correspond 
to the irreducible
components of the curve $C$.
\item[(b)] The edges correspond to the nodes.
\item[(c)] The markings
correspond to the marked points $p_i$.
\end{enumerate}
\item[(ii)] 
The genus function,  $g_{\tau}: V_{\tau} \rarr \mathbb{Z}^{\geq0}$,
where
$g_\tau(v)$ is the geometric genus of the corresponding
component of $C$. 
\item[(iii)]
The $H_2(X,\mathbb{Z})$-structure,
$\beta_{\tau}: V_\tau\rarr H_2(X, \beta)$, where
$\beta_{\tau}(v)$ equals
the $\mu$ push-forward
of the fundamental class 
of the corresponding component of $C$. 
\end{enumerate}

Following [BM], define 
$M_{\tau,n}(X,\beta_\tau)$ to
be the moduli space of
maps $\mu$ {\em together} with an isomorphism of
$\tau_\mu$  with a fixed 
stable graph $\tau$. The space $\overline{M}_{\tau,n}(X,\beta_\tau)$
is the compactification via stable maps where the 
vertices of $V_\tau$ may correspond to nodal curves.
Note $M_{\tau,n}(X,\beta_\tau)$ may not be dense in
$\overline{M}_{\tau,n}(X,\beta_\tau)$.

There is a canonical morphism
$$\pi_\tau: \overline{M}_{\tau,n}(X, \beta_\tau) \rarr
\overline{M}_{g,n}(X, \beta).$$
As $\tau$ varies over possible graphs,
the images of $\pi_\tau$ determine a (closed) decomposition
of the moduli space of maps.

Let $\tau$ be a stable, pointed, modular
graph with $H_2(\GG/\PP, \mathbb{Z})$-structure. 
The connectedness of 
$\overline{M}_{\tau, n}(\GG/\PP, \beta_\tau)$
will now be established.

\bigskip

\noindent {\em Proof of Theorem 2.} 
If $\beta_\tau=0$, the irreducibility of
$\overline{M}_{\tau, n}(\GG/\PP, \beta_\tau)$ is 
a direct 
consequence of the irreducibility of the
corresponding stratum in $\overline{M}_{g,n}$ and the
irreducibility of $\GG/\PP$. We may thus assume $\beta_\tau\neq 0$.

Fix the $\com^*$-action on $\GG/\PP$ as studied in Section \ref{cflow}.
Consider an arbitrary point
 $$[\mu] \in \overline{M}_{\tau,n}(\GG/\PP, \beta_\tau).$$
By the Kleiman-Bertini Theorem, a general $\GG$-translate $f$ of
$\mu$ satisfies conditions (i-iii) of Section \ref{cflow}.
As $\GG$ is connected, $[\mu]$ is connected to its general
$\GG$-translate $[f]$.

The point  $[f]$ is connected to
the limit:
$$[F] = {\text{lim}}_{t\rarr 0} [tf].$$
To prove the connectedness of $\overline
{M}_{\tau,n}(\GG/\PP, \beta_\tau)$,
it suffices to prove the set of limits $F$ lies in a
connected locus of the moduli space.
We will first construct the required connected locus of
$\overline
{M}_{\tau,n}(\GG/\PP, \beta_\tau)$.

The pair $(\tau, \beta_\tau)$ canonically determines a
family of maps $\gamma_b$ with nodal domains over 
a base $b\in B$.
For $v\in V_\tau$, let $\beta_\tau(v)= \sum_i a^v_i [P_i]$.
Define the base space $B$ as follows:
$$B = \prod_{v\in V_\tau} \overline{M}_{g(v),\text{val}
(v)+ \sum_i a_i^v},$$
where $\text{val}(v)$ is the valence of $v$ in $\tau$
(including nodes and markings). The extra $\sum_i a^v_i$ markings
each correspond to a basis homology element -- with
$a^v_j$ of these markings corresponding to $[P_j]$.
The degenerate cases $\overline{M}_{0,1}$ and $\overline{M}_{0,2}$
in the product $B$ are taken to be points.
$B$ is irreducible and hence connected.

For $b= \prod_v [b_v]\in B$, let
$$\gamma_b:D_b \rarr
\GG/\PP$$ 
be defined as follows:
\begin{enumerate}
\item[(i)] $D_b$ is obtained by attaching
the curves  $b_v$ by
connecting nodes as specified by $\tau$ and further
attaching $\proj^1$'s to each of the extra
points $\sum_i a^v_i$.
\item[(ii)] For each subcurve $b_v\subset D_b$, $\gamma_b(b_v)=s$.
\item[(iii)] For each $\proj^1$ corresponding
to $[P_j]$, $\gamma_b(\proj^1) \eqq P_j$.
\end{enumerate}
The family of maps $\gamma_b$ over $B$ then 
defines a morphism (via stabilization): 
$$\epsilon: B\rarr \overline{M}_{\tau,n}(\GG/\PP,\beta_\tau).$$
Certainly the image variety $\epsilon(B)$ is connected.

By Proposition \ref{pr2}, the limit $F$ is simply the stabilization
of $[\tilde{F}]$. Since $\tilde{F}= \gamma_b$ for some $b$,
the set of limits $F$ lies in a connected locus of 
$\overline{M}_{\tau,n}(\GG/\PP,\beta_\tau)$.
This concludes the proof of Theorem 2.
\epf

Theorem 1 is a special case of Theorem 2 (where $\tau$ has
a single vertex).
Corollary 2 is a simple consequences of Theorem 2.

\bigskip

\noindent{\em Proof of Corollary 2.}
In the genus 0 case, $\tau$ is a tree with genus function identically
zero. The moduli stack 
\begin{equation}
\label{rrrr}
\overline{\mathcal{M}}_{\tau,n}
(\GG/\PP, \beta_\tau)
\end{equation}
is constructed as a fiber product over the evaluation maps obtained
from the edges of $\tau$. 
We will prove $\overline{\mathcal{M}}_{\tau,n}(\GG/\PP, \beta_\tau)$
is a nonsingular Deligne-Mumford stack by induction on the
number of vertices of $\tau$.

First, suppose $\tau$ has only 1 vertex $v$. Then, the moduli
stack (\ref{rrrr}) is
$\overline{\mathcal{M}}_{0,\text{val}(v)}(\GG/\PP, \beta_\tau(v))$ 
-- a nonsingular
moduli stack by the convexity of $\GG/\PP$. 

Next, let $\tau$ have $m$ vertices and let $v$ be an extremal
vertex ($v$ is incident to exactly 1 edge).
Let
$p\in \GG/\PP$ be a point. By the Kleiman-Bertini Theorem,
\begin{equation}
\label{ddd}
\text{ev}_1^{-1}(p) \subset \overline{\mathcal{M}}_{0,\text{val}(v)}
(\GG/\PP, \beta_\tau(v))
\end{equation}
is a nonsingular Deligne-Mumford stack for the general point $p$ (and hence
{\em every} point $p$).
Let $\tau'$ be the graph obtained by removing $v$ from $\tau$
and adding an extra marking corresponding to the broken node.
The moduli stack (\ref{rrrr}) is fibered over
\begin{equation}
\label{vbnn} 
\overline{\mathcal{M}}_{\tau', n'+1}(\GG/\PP, \beta_\tau')
\end{equation}
with fiber (\ref{ddd}). 
As (\ref{vbnn}) is nonsingular by induction, the stack (\ref{rrrr})
is thus nonsingular. This completes the induction step.

Finally, since $\overline{\mathcal{M}}_{\tau,n}(\GG/\PP, \beta_\tau)$
is a nonsingular and connected Deligne-Mumford stack, it is
irreducible.
\epf

\section{Rationality}
We first review a basic rationality result proven in
a sequence papers by Katsylo and Bogomolov [Ka], [Bog].
Let $V= \com^2$ be a vector space.
Let $a_1,a_2, \ldots, a_n$ be a sequence of positive integers
with $\sum_i a_i \geq 3$. Then, the quotient
\begin{equation}
\label{dfgdfgdfg}
\proj(\text{Sym}^{a_1} V^*) \times \cdots \times
\proj(\text{Sym}^{a_n} V^*) \ // \ \mathbf{PGL}(V)
\end{equation}
is a rational variety -- we may take any non-empty invariant theory
quotient. Geometrically, the quotient
(\ref{dfgdfgdfg}) is birational to the moduli space
quotient
\begin{equation}
M_{0,\sum_i a_i}\ / \ \Sigma_{a_1} \times \Sigma_{a_2} \times
\cdots \times \Sigma_{a_n}
\end{equation}
where $\Sigma$ is the symmetric group.
Essentially, the rationality of (\ref{dfgdfgdfg}) is deduced
from rationality in case $n=1$ [Ka]. Proofs in the
$n=1$ case may be found in [Ka], [Bog].

We will also need the following simple Lemma.
\begin{lm}
\label{qwwq}
Let $W$ be any finite dimensional linear representation of
$\mathbf{A}$ where $\mathbf{A}= \Sigma_2$ or $\mathbf{A}=\Sigma_3$.
Then, $W/\mathbf{A}$ is rational.
\end{lm}
\bpf
By the complete reducibility of representations and the
fact that a $\mathbf{GL}$-bundle is locally
trivial in the Zariski topology, it suffices to prove
the Lemma in case $W$ is an irreducible representation.
It is then easily checked by hand the two irreducible representation
of $\Sigma_2$ and the three
 irreducible representations of $\Sigma_3$
have rational quotients.
\epf

\bigskip 

\noindent{\em Proof of Theorem 3.}
Fix the $\com^*$-action on $\GG/\PP$ as studied in Section \ref{cflow}.
We first consider the moduli space
$$\overline{M}=\overline{M}_{0,n}(\GG/ \PP, \beta=\sum_i a_i [P_i])$$
where the property 
\begin{equation}
\label{dga}
n+ \sum_i a_i \geq 4
\end{equation}
is satisfied.

Let $\tau$ be the graph with a single vertex $v$ with $n$
markings, and let $\beta_\tau(v)= \sum_i a_i[P_i]$.
Let $\gamma_b$ over $B$ be the family of maps constructed
canonically from $(\tau, \beta_\tau)$ in the proof of
Theorem 2.
The base $B$ is simply:
\begin{equation}
\label{ggff}
B= \overline{M}_{0,n+\sum_i a_i}.
\end{equation}
The map $\gamma_b$ over a general point $b\in B$ has no
map automorphisms (as $n+\sum_i a_i \geq 4$). 
Hence, the image $\epsilon(B)$ in $\overline{M}$ intersects
the nonsingular (automorphism-free) locus of the moduli space
$\overline{M}^0\subset \overline{M}$. Let
$$\epsilon(B)^0= \epsilon(B) \cap \overline{M}^0,$$
and let $B^0= \epsilon^{-1}(\epsilon(B)^0).$
The map
$$B^0 \rarr \epsilon(B)^0$$
is simply a quotient of $B^0$ by the
natural $\Sigma_{a_1} \times \cdots \times
\Sigma_{a_n}$ action on (\ref{ggff}).
By the rationality result (\ref{dfgdfgdfg}), $\epsilon(B)^0$
is rational.

Consider now the $\com^*$-action 
on $\overline{M}^0$ by translation.
As $\overline{M}^0$ is a nonsingular, irreducible, quasi-projective
variety, we may study the Bialynicki-Birula stratification
of $\overline{M}^0$. By the proof of Theorem 2, 
$\epsilon(B)^0$ is a
$\com^*$-fixed locus which contains the
limit, 
$$\text{lim}_{t\rarr 0} \ t[f],$$ of the
general point $[f]\in \overline{M}^0$. By [Bi], $\overline{M}^0$
is
birational to an affine bundle over $\epsilon(B)^0$.
Therefore, $\overline{M}$ is rational.
The proof of Theorem 3 is  complete in case (\ref{dga})
is satisfied.

Next, we will consider the case where the sum (\ref{dga}) is at most
3. In this case, the base $B$ is a point.
If $\epsilon(B)$ lies in the automorphism-free locus, the
previous argument proving the rationality
of $\overline{M}_{0,n}(\GG/\PP,\beta)$
is still valid. There are exactly four cases in
which the point $\epsilon(B)$ corresponds to a map with nontrivial
automorphisms:
\begin{enumerate}
\item[(i)] $n=0$, $\beta= 3[P_i]$.
\item[(ii)] $n=0$, $\beta= 2[P_i] + [P_j]$, $i\neq j$.
\item[(iii)] $n=0$, $\beta=2[P_i]$.
\item[(iv)] $n=1$, $\beta= 2[P_i]$.
\end{enumerate}
Here, the Deligne-Mumford
stack structure of these moduli spaces is important.
The automorphism group in case (i) is $\Sigma_3$
and in cases (ii-iv) is $\Sigma_2$.
In each case, we will show the coarse
moduli space $\overline{M}_{0,n}(\GG/\PP,
\beta)$ 
is birational to a quotient of a linear representation of the
corresponding automorphism group.

Consider first the case (i): $n=0$, $\beta=3[P_i]$.
Let $\epsilon(B)=[\gamma]$. Let $[\mu]$ denote
the unique 3-pointed stable map obtained from
$\gamma$ by marking each $\proj^1 \eqq P_i$ by
a point lying over $x_i$. Certainly,
$[\mu]\in \overline{M}^0_{0,3}(\GG/\PP, \beta)$

We will study: 
$$N\subset \overline{M}^0_{0,3}(\GG/\PP, \beta)$$
where $N$ is the component of the locus of transverse intersection
of the three divisors $\text{ev}_1^{-1}(\overline{A}_i)$, 
$\text{ev}_2^{-1}(\overline{A}_i)$,
and $\text{ev}_3^{-1}(\overline{A}_i)$ containing $[\mu]$.
The torus $\com^*$ acts on $N$ by translation. 
By an argument exactly parallel to the flow result of
Proposition \ref{pr2}, we deduce 
$$\text{lim}_{t\rarr 0} t [f] = [\mu]$$
for a general element $[f]\in N$. As $N$ is a nonsingular, quasi-projective
scheme, Theorem 2.5 of [Bi] implies that
$N$ is $\com^*$-equivariantly birational to the
tangent $\com^*$-representation at $[\mu]$.

There is a 
$\Sigma_3$-action on $N$ by permutation of the markings.
The $\com^*$ and $\Sigma_3$ actions commute.
A slightly refined version of Theorem 2.5 of [Bi]
shows $N$ is $\com^* \times \Sigma_3$-equivariantly birational
to the tangent $\com^* \times \Sigma_3$-representation
at $[\mu]$. Lemma \ref{blw} below explains the refinements of
the results of [Bi] needed here.
$N/\Sigma_3$ is birational to $\overline{M}_{0,0}(\GG/\PP, \beta)$.
Hence, by Lemma \ref{qwwq}, Theorem 3 is proven in case (i).

A similar strategy is used in cases (ii-iv).
In each of these cases, let 
$\epsilon(B)=[\gamma]$ and let $[\mu]$ denote the
rigidification by adding 2 new markings ${\bullet, \bullet'}$ which lie over
$x_i$. 
The locus $N$ is chosen 
as the corresponding transverse intersection locus
of $\text{ev}_{\bullet}^{-1}(\overline{A}_i)$ 
and $\text{ev}_{\bullet'}^{-2}(\overline{A}_i)$ in
the maps space with the new markings.
$N$ is then $\com^*\times \Sigma_2$-equivariantly
birational to the tangent $\com^* \times \Sigma_2$-representation
of $N$ at $[\mu]$ by the refined Lemma \ref{blw}. Theorem 3 is then
a consequence of Lemma \ref{qwwq} since $N/\Sigma_2$
is birational to the moduli space of maps considered in the case.
\epf

\begin{lm}
\label{blw}
Let $\mathbf{A}$ be a finite group.
Let $S$ be a nonsingular, irreducible,
quasi-projective scheme with
a $\com^* \times \mathbf{A}$-action and a $\com^* \times \mathbf{A}$-fixed
point $s\in S$. Let $T_s$ denote the $\com^* \times \mathbf{A}$-representation
on the tangent space at $s$.
Suppose the  $\com^*$-action is
fully definite at $s$. Then, there is 
$\com^*\times \mathbf{A}$-equivariant isomorphism
between an open set of  $(S,s)$ and $(T_s,0)$.
\end{lm}

\bpf
We note $\com^* \times \mathbf{A}$ is a linearly reductive group.
By Theorem 2.4 of [Bi] for linearly reductive group actions, we may find
a third nonsingular irreducible
pointed space $(Z,z)$ with a $\com^* \times \mathbf{A}$-action
and equivariant, \'etale, morphisms:
$$\pi_1: (Z,z) \rarr (S,s),$$
$$\pi_2:(Z,z) \rarr(T_s,s).$$
In the proof of Theorem 2.5 of [Bi], such morphisms $\pi_1$ and $\pi_2$
are proven to be open immersions {\em by a study of only the
$\com^*$-action}. Hence, the morphisms $\pi_1$ and $\pi_2$ are
open immersions in our case.
By the full definiteness of the $\com^*$-representation
on $T_s$, the morphism $\pi_2$ is then an isomorphism. 
\epf

\end{document}